\documentclass[12pt]{article}
\usepackage{amsmath}

\usepackage[margin=1in, centering]{geometry}

\usepackage{placeins}
\usepackage{verbatim}

\usepackage{authblk}
\usepackage{url}
\usepackage{amsthm}

\usepackage{float}

\usepackage[dvips]{graphicx}
\usepackage{color}
\usepackage{amsfonts,amssymb,amsmath,}
\usepackage{mathrsfs}
\usepackage[T1]{fontenc}
\usepackage{ae}

\newtheorem{theorem}{Theorem}

\newtheorem{lem}{Lemma}

\numberwithin{equation}{section}
\numberwithin{table}{section}
\numberwithin{figure}{section}

\renewcommand{\(}{\left(}
\renewcommand{\)}{\right)}

\allowdisplaybreaks[4]

\begin{document}
\title{Summatory function of the number of prime factors}
\author{Xianchang Meng}
\date{}

\maketitle
\begin{abstract}
 We consider the summatory function of the number of prime factors for integers $\leq x$ over  arithmetic progressions. Numerical experiments suggest that  some arithmetic progressions consist more number of prime factors than others. Greg Martin conjectured that the difference of the summatory functions should attain a constant sign for all sufficiently large $x$.  In this paper, we provide strong evidence for Greg Martin's conjecture. Moreover, we derive a general theorem for arithmetic functions from the Selberg class.

\end{abstract}

\let\thefootnote\relax

\let\thefootnote\relax

\section{Introduction and statement of results}
Let $\Omega(n)$ be the number of prime factors of $n$ counted with multiplicity, and $\omega(n)$ be the number of distinct prime factors of $n$. Hardy and Ramanujan \cite{Hardy-Rama} showed that
\begin{equation}\label{sum-l-omega}
\sum_{n\leq x}\omega(n)=x\log\log x+Ax+O\(\frac{x}{\log\log x} \), 
\end{equation}
and 
\begin{equation}\label{sum-B-Omega}
\sum_{n\leq x}\Omega(n)=x\log\log x+Bx+O\(\frac{x}{\log\log x} \), 
\end{equation}
for some constants $A$ and $B$. Thus, the average order of $\Omega(n)$ and $\omega(n)$ is $\log\log n$. The famous Erd\H os-Kac Theorem  (see \cite{Erdos-Kac} or \cite{Renyi-Turan}) says that the limiting distribution of 
$$ \frac{\omega(n)-\log\log n}{\sqrt{\log\log n}}$$
is the standard normal distribution. 

One may get better error terms for the summatory functions in \eqref{sum-l-omega} and \eqref{sum-B-Omega} by assuming the Riemann Hypothesis. Wolke \cite{Wolke} proved that the Riemann Hypothesis is true if and only if 
\begin{equation}\label{wolke-RH}
\sum_{n\leq x} \omega(n)=x\sum_{0\leq j\leq \log x/2}\frac{P_j(\log\log x)}{\log^j x}+O\( x^{\frac{1}{2}+\epsilon}\), 
\end{equation}
where $P_j$ are some polynomials of degree $1$.  Similarly, under the Extended Riemann Hypothesis (ERH) for Dirichlet $L$-functions, we expect the following formula for the corresponding summatory function over arithmetic progressions,
$$\sum_{\substack{n\equiv a \bmod q\\n\leq x}} \omega(n)=\frac{x}{\phi(q)} \sum_{0\leq j\leq \log x/2}\frac{P_j(\log\log x)}{\log^j x}+O\( x^{\frac{1}{2}+\epsilon}\), \quad (a, q)=1.$$
Thus, when we compare two different arithmetic progressions $a\not\equiv b \bmod q$, $(a, q)=1$ and $(b, q)=1$, under the ERH, we have
$$ \sum_{\substack{n\equiv a \bmod q\\n\leq x}} \omega(n)-\sum_{\substack{n\equiv b \bmod q\\n\leq x}} \omega(n)=O\(x^{\frac{1}{2}+\epsilon} \).$$
Based on the above analysis, one may expect that there are roughly equal number of prime factors in different arithmetic progressions. 

However, Greg Martin noticed that some arithmetic progressions actually  have more number of prime factors by doing numerical calculations. He conjectured that
\begin{equation}\label{Martin-Conj}
 \sum_{\substack{n\equiv 1 \bmod 4\\n\leq x}} \omega(n)-\sum_{\substack{n\equiv 3 \bmod 4\\n\leq x}} \omega(n)<0
\end{equation}
holds for all sufficiently large $x$. 

In this paper, we give strong evidence for Greg Martin's conjecture \eqref{Martin-Conj}.  We consider both cases of $\Omega(n)$ and $\omega(n)$. Let $\chi\neq \chi_0$ be a non-principal Dirichlet character modulo $q$, and denote
\begin{equation*}
\psi_{f}(x,\chi):=\sum_{n\leq x}\chi(n)f(n),
\end{equation*}
where $f=\omega$ or $\Omega$.

\begin{theorem}\label{thm-omega}
	Assume the ERH, the zeros of $L(s, \chi)$ are simple, and $L(\frac{1}{2}, \chi)\neq 0$. For any fixed large $T_0$, we have
	\begin{align}
	\psi_{\omega}(x, \chi)=&-a(\chi)\left\{L\(\frac{1}{2}, \chi\) \frac{\sqrt{x}}{\log x} +\(2L\(\frac{1}{2}, \chi\)-L'\(\frac{1}{2}, \chi\)\) \frac{\sqrt{x}}{\log^2 x}\right\}\nonumber\\
	&+\frac{\sqrt{x}}{\log^2 x}\left\{\sum_{|\gamma|\leq T_0} \frac{L'(\rho, \chi) x^{i\gamma}}{\frac{1}{2}+i\gamma}
	+\Sigma(x, T_0)\right\},
	\end{align}	
	where $a(\chi)=1$ if $\chi$ is real, 0 otherwise, and $\forall \epsilon >0$,
	\begin{equation}
	\limsup_{Y\rightarrow\infty}\frac{1}{Y}\int_1^Y \left| \Sigma(e^y, T_0) \right|^2 dy\ll \frac{1}{T_0^{1-\epsilon}}.
	\end{equation}
\end{theorem}

Similarly, we have the following result for $\psi_{\Omega}(x, \chi)$. 
\begin{theorem}\label{thm-Omega}
	Assume the ERH, the simplicity of zeros, and $L(\frac{1}{2}, \chi)\neq 0$. For any fixed large $T_0$, 
	\begin{align}
	\psi_{\Omega}(x, \chi)=&a(\chi)\left\{L\(\frac{1}{2}, \chi\) \frac{\sqrt{x}}{\log x} +\(2L\(\frac{1}{2}, \chi\)-L'\(\frac{1}{2}, \chi\)\) \frac{\sqrt{x}}{\log^2 x}\right\}\nonumber\\
	&+\frac{\sqrt{x}}{\log^2 x}\left\{\sum_{|\gamma|\leq T_0} \frac{L'(\rho, \chi) x^{i\gamma}}{\frac{1}{2}+i\gamma}
	+\Sigma(x, T_0)\right\},
	\end{align}	
	where $a(\chi)=1$ if $\chi$ is real, 0 otherwise, and $\forall \epsilon >0$,
	\begin{equation}
	\limsup_{Y\rightarrow\infty}\frac{1}{Y}\int_1^Y \left| \Sigma(e^y, T_0) \right|^2 dy\ll \frac{1}{T_0^{1-\epsilon}}.
	\end{equation}
\end{theorem}

The Linear Independence Conjecture (LI) states that all the positive imaginary parts of zeros of $L(s, \chi)$ are linearly independent over $\mathbb{Q}$. In \cite{Meng-kfree} and \cite{Meng}, the author used ERH and LI to study the the distribution of $k$-free numbers and subtle difference of the number of products of $k$ primes among different arithmetic progressions. 

Let $\chi_{-4}$ be the non-principal Dirichlet character $\bmod 4$, i.e. $\chi_{-4}(1)=1$ and  $\chi_{-4}(3)=-1$. Let 
$$P_{\omega}:= \Bigg\{   N\in\mathbb{N}:   \sum_{\substack{n\equiv 1 \bmod 4\\n\leq N}} \omega(n)-\sum_{\substack{n\equiv 3 \bmod 4\\n\leq N}} \omega(n)<0   \Bigg\},$$
and 
$$P_{\Omega}:= \Bigg\{   N\in\mathbb{N}:   \sum_{\substack{n\equiv 1 \bmod 4\\n\leq N}} \Omega(n)-\sum_{\substack{n\equiv 3 \bmod 4\\n\leq N}} \Omega(n)>0   \Bigg\}.$$
Since $L\(\frac{1}{2}, \chi_{-4}\)$>0, by Theorems \ref{thm-omega} and \ref{thm-Omega} and further assuming LI, we deduce that the \textit{logarithmic densities} of the sets $P_{\omega}$ and $P_{\Omega}$ are both equal to $1$. The \textit{logarithmic density} of a set $S$ is defined to be
$$\delta(S):=\lim_{X\rightarrow\infty} \frac{1}{\log X} \sum_{n\leq X, n\in S} \frac{1}{n}.$$

Although we cannot prove the full Greg Martin's conjecture  (\ref{Martin-Conj}), we give strong evidence to support this conjecture subject to the ERH and LI. In order to prove the full conjecture, one may need to formulate new ideas and introduce more powerful tools. 

\begin{figure}[H]
	\centering
	\includegraphics[width=0.8\linewidth]{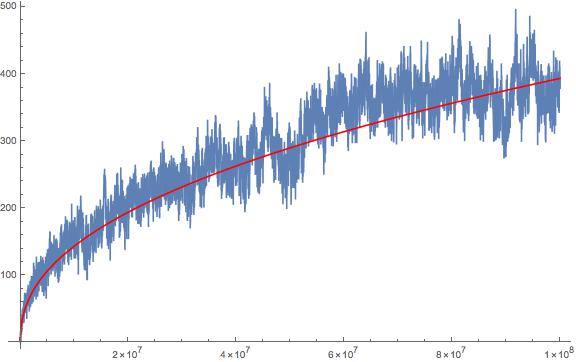}\\
	\caption{Graph of $\psi_{\Omega}(x, \chi_{-4})$}\label{fig-Omega}
\end{figure}

\begin{figure}[H]
	\centering
	\includegraphics[width=0.8\linewidth]{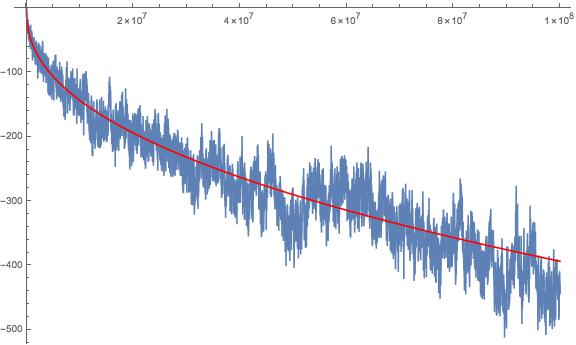}\\
	\caption{Graph of $\psi_{\omega}(x, \chi_{-4})$}\label{fig-omega}
\end{figure}

We numerically calculated $\psi_{\Omega}(x, \chi_{-4})$ and $\psi_{\omega}(x, \chi_{-4})$ for $x\leq 10^8$. See Figures \ref{fig-omega} and \ref{fig-Omega}. In these two graphs, the blue dots are numerical data and the red lines are the functions of the main terms predicted in our theorems. We see that our theorems match very well with  numerical calculations. 

\begin{figure}[H]
\centering
\includegraphics[width=0.8\linewidth]{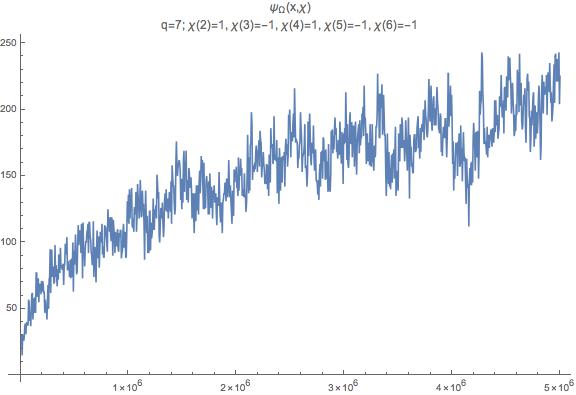}

\includegraphics[width=0.8\linewidth]{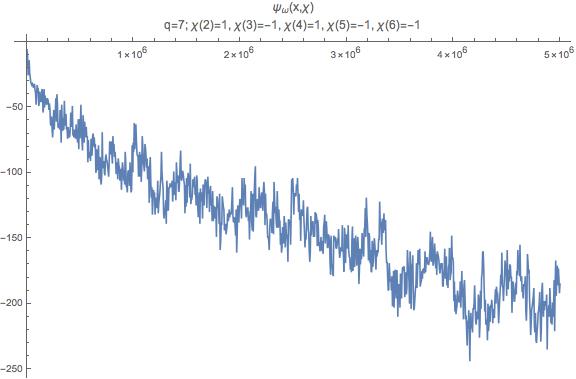}
\caption{Real Dirichlet character mod $7$}\label{real-character}

\end{figure}

\begin{figure}[H]
\centering
\includegraphics[width=0.8\linewidth]{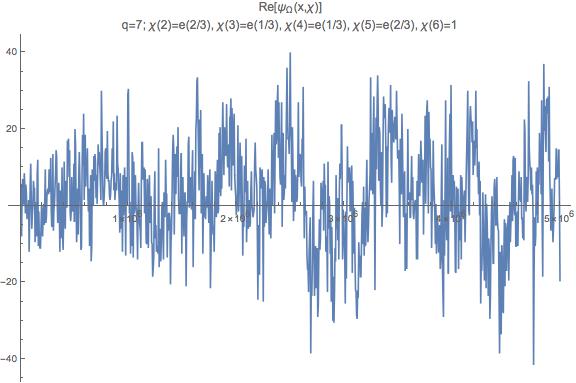}
\includegraphics[width=0.8\linewidth]{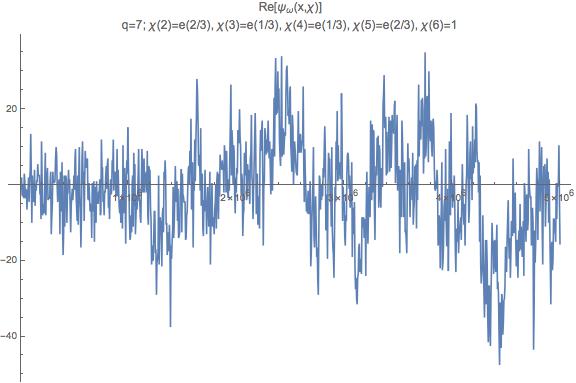}
\caption{Complex Dirichlet character mod $7$}\label{complex-character}

\end{figure}

Greg Martin provided us more numerical data for different Dirichlet characters, including both real and complex characters. See Figures \ref{real-character} and \ref{complex-character}.

The structure of this paper is as follows. In Section \ref{Sec-proof}, we give the proof of our main theorems. In Section \ref{Sec-Selberg}, we generalize the method used in \cite{Meng} and this paper to more general arithmetic functions in Selberg class.

\section{Proof of Theorems}\label{Sec-proof}

Consider the Dirichlet series,
\begin{equation}\label{F-omega}
F_{\omega}(s, \chi):=\sum_{n=1}^{\infty} \frac{\omega(n)\chi(n)}{n^s}=L(s, \chi)F(s, \chi),
\end{equation}
and 
\begin{equation}\label{F-Omega}
F_{\Omega}(s, \chi):=\sum_{n=1}^{\infty}\frac{\Omega(n)\chi(n)}{n^s}=L(s, \chi)\(F(s, \chi)+F(2s, \chi^2)+F(3s, \chi^3)+\cdots \),
\end{equation}
where
\begin{equation}
F(s, \chi)=\sum_{p~\text{prime}} \frac{\chi(p)}{p^s}=\log L(s, \chi)-\frac{1}{2}\log L(2s, \chi^2)+G(s),
\end{equation}
with $G(s)$ absolutely convergent for $\Re(s)\geq \sigma_0=0.34$.

In our proof, we borrow the idea developed in \cite{Meng} to deal with the singularities of $F_{\omega}(s, \chi)$ and $F_{\Omega}(s, \chi)$. In \cite{Meng}, we mainly dealt with powers of $\log L(s, \chi)$. In this paper, we have to take care of $L(s, \chi)\log L(s, \chi)$ with more careful calculations.

\subsection{Contour Integral Representation}
By Perron's formula  (\cite{Kara}, Chapter V, Theorem 1), we have
\begin{lem}\label{lem-psi-Perron}
	For any $T\geq 2$,
	\begin{equation*}
	\psi_{f}(x,\chi)=\frac{1}{2\pi i}\int_{c-iT}^{c+iT} F_{f}(s, \chi)\frac{x^s}{s} ds +O\(\frac{x\log x}{T}+1\),
	\end{equation*}
	where $c=1+\frac{1}{\log x}$, and $f=\omega$ or $\Omega$.
\end{lem}

Under the ERH, using the similar method as in \cite{Titch} (Theorem 14.16), one can show that, for any $\epsilon>0$ and $\forall \chi\neq \chi_0 \bmod q$, there exists a sequence of numbers $\mathcal{T}=\{T_n \}_{n=0}^{\infty}$ satisfying $n\leq T_n\leq n+1$ such that,
$T_n^{-\epsilon}\ll |L(\sigma+iT_n, \chi)|\ll T_n^{\delta+\epsilon}, ~(\frac{1}{2}-\delta<\sigma<2).$

Let $\rho$ be a zero of $L(s, \chi)$, $\Delta_{\rho}$ be the distance of $\rho$ to the nearest other zero, and $D_{\gamma}:=\min\limits_{T\in \mathcal{T}}(|\gamma-T|)$.
For each zero $\rho$, and $X>0$, let $\mathcal{H}(\rho, X)$ denote the truncated Hankel contour surrounding the point $s=\rho$ with radius $0<r_{\rho}\leq \min(\frac{1}{x}, \frac{\Delta_{\rho}}{3}, \frac{D_{\gamma}}{2})$, which includes the circle $|s-\rho|=r_{\rho}$ excluding the point $s=\rho-r_\rho$, and the half-line $(\rho-X, \rho-r]$ traced twice with arguments $+\pi$ and $-\pi$ respectively. Let $\mathcal{H}(\frac{1}{2}, X)$ denote the corresponding Hankel contour surrounding $s=\frac{1}{2}$ with radius $r_0=\frac{1}{x}$.

\begin{figure}[h!]
	\centering
	\includegraphics[width=10cm]{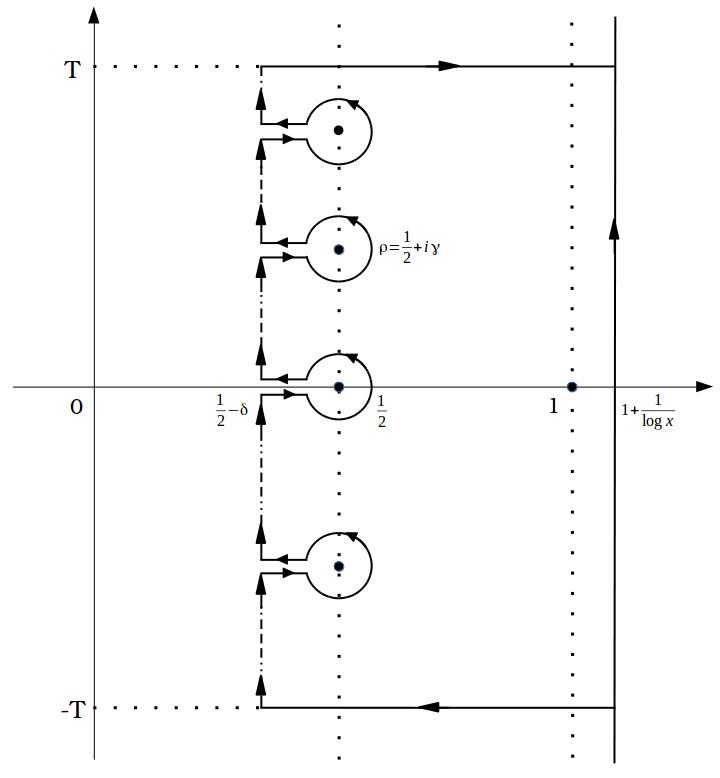}\\
	\caption{Integration contour}\label{contour}
\end{figure}

Take $\delta=\frac{1}{100}$. By Lemma \ref{lem-psi-Perron}, we pull the contour to the left to the line $\Re(s)=\frac{1}{2}-\delta$ using the truncated Hankel contour  $\mathcal{H}(\rho, \delta)$ to avoid the zeros of $L(s, \chi)$ and using $\mathcal{H}(\frac{1}{2}, \delta)$ to avoid the point $s=\frac{1}{2}$. See Figure \ref{contour}.

Then we have the following lemma. 
\begin{lem}\label{lem-psi-contour}
	Assume the ERH, and $L(\frac{1}{2}, \chi)\neq 0$ $(\chi\neq \chi_0)$. Then for $T\in \mathcal{T}$, 
	\begin{align}
	\psi_f(x, \chi)=&\sum_{|\gamma|\leq T} \frac{1}{2\pi i} \int_{\mathcal{H}(\rho, \delta)} F_f(s, \chi)\frac{x^s}{s}ds +a(\chi) \frac{1}{2\pi i} \int_{\mathcal{H}(\frac{1}{2}, \delta)} F_f(s, \chi)\frac{x^s}{s}ds\nonumber\\
	&+O\(\frac{x\log x}{T}+\frac{x}{T^{1-\delta-\epsilon}}+x^{\frac{1}{2}-\delta}T^{\delta+\epsilon} \),
	\end{align}
	where $a(\chi)=1$ if $\chi$ is real, $0$ otherwise, and $f=\omega$ or $\Omega$. 
\end{lem}
\noindent\textbf{\textit{Proof.}} By (\ref{F-omega}) and (\ref{F-Omega}), if $\chi$ is not real, $s=\frac{1}{2}$ is not a singularity. Under the ERH, the integral on the horizontal line is bounded by 
\begin{equation}
O\(T^{\delta+\epsilon}\int_{\frac{1}{2}-\delta}^{c} \frac{x^{\sigma}}{T}d\sigma \)=O\(\frac{x}{T^{1-\delta-\epsilon}} \).
\end{equation}
And the integral on the vertical line $\Re (s)=\frac{1}{2}-\delta$ is bounded by
\begin{equation}
O\(\int_{-T}^T \frac{(|t|+2)^{\delta+\epsilon} x^{\frac{1}{2}-\delta}}{|t|+2}dt \)=O\(x^{\frac{1}{2}-\delta}T^{\delta+\epsilon} \).
\end{equation}
Then by Lemma \ref{lem-psi-Perron}, we get the desired formula. \qed

\vspace{1em} 
Let $\mathcal{H}(0, X)$ be the truncated Hankel contour surrounding $0$ with radius $r$. For simplicity,  we denote
$$\frac{1}{\Gamma_j(u)}:=\left[\frac{d^j}{dz^j}\(\frac{1}{\Gamma(z)}\)\right]_{z=u}.$$
\begin{lem}[\cite{Lau-Wu}, Lemma 5]\label{lem-Gamma-j} For $X>1$, $z\in \mathbb{C}$ and $j\in \mathbb{Z}^{+}$, we have
	\begin{equation*}
	\frac{1}{2\pi i}\int_{\mathcal{H}(0, X)} w^{-z}(\log w)^j e^w dw=(-1)^j \frac{d^j}{dz^j}\(\frac{1}{\Gamma(z)}\)  +E_{j, z}(X),
	\end{equation*}
	where
	\begin{equation*}
	|E_{j,z}(X)|\leq \frac{e^{\pi|\Im(z)|}}{2\pi}\int_{X}^{\infty} \frac{(\log t+\pi)^j}{t^{\Re(z)}e^t}dt.
	\end{equation*}
\end{lem}

\begin{lem}[\cite{Meng}, Lemma 12]\label{lem-Meng-8}
	For any integers $k\geq 1$ and $m\geq 0$. We have
	\begin{equation*}
	\int_0^{\delta} |(\log \sigma -i\pi)^k-(\log\sigma +i\pi)^k|\sigma^m x^{-\sigma}d\sigma\ll_{m,k} \frac{(\log\log x)^{k-1}}{(\log x)^{m+1}}.
	\end{equation*}
\end{lem}

Similar to the proof of Lemma 8 in \cite{Meng}, we have the following result.
\begin{lem}\label{lem-main-term}
	Let $\mathcal{H}(a, \delta)$  be the truncated Hankel contour surrounding a complex number $a \ (\Re(a)>2\delta)$ with radius $0<r\ll \frac{1}{x}$. Then, for any integers $k\geq 1$, $l\geq 0$
	\begin{align*}
	&\frac{1}{2\pi i}\int_{\mathcal{H}(a, \delta)} (s-a)^l \log^k(s-a)\frac{x^s}{s}ds\nonumber\\
	&= \frac{(-1)^k x^a}{a(\log x)^{l+1}}  \sum_{j=1}^k {k \choose j} (\log\log x)^{k-j} \frac{1}{\Gamma_j(-l)}-\frac{(-1)^k x^a}{a^2(\log x)^{l+2}} \sum_{j=1}^k {k \choose j} (\log\log x)^{k-j} \frac{1}{\Gamma_j(-l-1)}\nonumber\\
	&\quad +O_{k,l}\( \frac{|x^a|}{|a|^2|\Re(a)-\delta|}\frac{(\log\log x)^{k-1}}{(\log x)^{l+3}}\)+O_{k, l}\(\frac{|x^{a-\delta/2}|}{|a|}+\frac{|x^{a-\delta/2}|}{|a|^2}\).
	\end{align*}
\end{lem}
\noindent\textbf{\textit{Proof of Lemma \ref{lem-main-term}.}} We have the equality
\begin{equation*}
\frac{1}{s}=\frac{1}{a}+\frac{a-s}{a^2}+\frac{(a-s)^2}{a^2 s}.
\end{equation*}
With the above equality, we write the integral in the lemma as
\begin{equation*}
\frac{1}{2\pi i}\int_{\mathcal{H}(a, \delta)} (s-a)^l\log^k(s-a)\(\frac{1}{a}+\frac{a-s}{a^2}+\frac{(a-s)^2}{a^2 s}\)x^s ds =: I_1+I_2+I_3.
\end{equation*}
For $I_3$, using Lemma \ref{lem-Meng-8}, we get
\begin{align}\label{pf-lem-main-I3}
&\int_{\mathcal{H}(a,\delta)} (s-a)^l\log^k(s-a)\frac{(s-a)^{2}}{a^2 s} x^s ds\nonumber\\
&\leq \left|\int_r^{\delta}\left((\log \sigma-i\pi)^k-(\log\sigma+i\pi)^k\right)\sigma^{l+2} x^{-\sigma}\frac{x^a}{a^2 (a-\sigma)}d\sigma\right|\nonumber\\
&\quad +\int_{-\pi}^{\pi} x^{\Re(a)+r} \(\log \frac{1}{r}+\pi\)^k \frac{r^{l+2}}{|a|^2 |\Re(a)-r|}  r d\alpha\nonumber\\
&\ll\frac{|x^a|}{|a|^2 |\Re(a)-\delta|} \(\int_0^{\delta} |(\log \sigma -i\pi)^k-(\log\sigma +i\pi)^k|\sigma^{l+2} x^{-\sigma}d\sigma+ \frac{(\log\frac{1}{r}+\pi)^{k}}{(1/r)^{l+3}}\)\nonumber\\
&\ll_{k, l} \frac{|x^a|}{|a|^2|\Re(a)-\delta|}\(\frac{(\log\log x)^{k-1}}{(\log x)^{l+3}}+ \frac{1}{x^{l+3-\epsilon}}\)\ll_{k, l} \frac{|x^a|}{|a|^2|\Re(a)-\delta|}\frac{(\log\log x)^{k-1}}{(\log x)^{l+3}}
\end{align}
For $I_1$, using change of variable $(s-a)\log x=w$, by Lemma \ref{lem-Gamma-j},
\begin{align}\label{pf-lem-main-I1}
I_1&=\frac{1}{2\pi i}\frac{1}{(\log x)^{l+1}}\int_{\mathcal{H}(0,\delta\log x)}w^l(\log w-\log\log x)^k \frac{x^a e^w}{a}dw\nonumber\\
&=\frac{x^a}{a(\log x)^{l+1}}  (-1)^k (\log\log x)^k \frac{1}{2\pi i}\int_{\mathcal{H}(0,\delta\log x)} w^l e^w dw \nonumber\\
& \quad +\frac{x^a}{a(\log x)^{l+1}}\sum_{j=1}^k {k\choose j} \frac{1}{2\pi i}  \int_{\mathcal{H}(0, \delta\log x)} (-\log\log x)^{k-j} w^l (\log w)^j e^w dw  \nonumber\\
&=\frac{(-1)^k x^a}{a(\log x)^{l+1}} \sum_{j=1}^k {k \choose j} (\log\log x)^{k-j} \frac{1}{\Gamma_j(-l)}  \nonumber \\
&  \quad +\frac{x^a}{a(\log x)^{l+1}} \sum_{j=1}^k {k \choose j} E_{j,-l}(\delta\log x) (-\log\log x)^{k-j}.
\end{align}
By Lemma \ref{lem-Gamma-j},
\begin{equation*}
|E_{j,-l}(\delta\log x)|\leq \frac{1}{2\pi}\int_{\delta\log x}^{\infty} \frac{t^l(\log t+\pi)^j}{e^t}dt
\ll_{j} e^{-\frac{\delta\log x}{2}} \int_{\frac{\delta\log x}{2}}^{\infty} \frac{t^l(\log t)^j}{e^{t/2}}dt
\ll_{j, l} x^{-\frac{\delta}{2}}.
\end{equation*}
Hence, we get
\begin{align}\label{pf-lem-main-E}
&\left| \frac{x^a}{a(\log x)^{l+1}} \sum_{j=1}^k {k \choose j} E_{j,-l}(\delta\log x) (-\log\log x)^{k-j} \right| \nonumber\\
&\quad \ll_{k,l} \frac{x^{\Re(a)}}{|a|(\log x)^{l+1}} \sum_{j=1}^k x^{-\frac{\delta}{2}} (\log\log x)^{k-j}\ll_{k,l} \frac{|x^{a-\delta/2|}}{|a|}.
\end{align}
For $I_2$, we have
\begin{equation}
\frac{1}{2\pi i}\int_{\mathcal{H}(a,\delta)} (s-a)^l\log^k(s-a)\frac{a-s}{a^2} x^s ds=-\frac{1}{a}\frac{1}{2\pi i}\int_{\mathcal{H}(a,\delta)} (s-a)^{l+1}\log^k(s-a) \frac{x^s}{a} ds.
\end{equation}
Then, we replace $l$ by $l+1$ in the formula of $I_1$ and multiply by $-\frac{1}{a}$. Thus, we get
\begin{align}\label{pf-lem-main-I2}
I_2=-\frac{(-1)^k x^a}{a^2(\log x)^{l+2}} \sum_{j=1}^k {k \choose j} (\log\log x)^{k-j} \frac{1}{\Gamma_j(-l-1)}+O_{k,l}\(\frac{|x^{a-\delta/2}|}{|a|^2} \).
\end{align}

Combining (\ref{pf-lem-main-I3}), (\ref{pf-lem-main-I1}), (\ref{pf-lem-main-I2}), and (\ref{pf-lem-main-E}), we get the conclusion of this lemma. \qed

\subsection{Proof of Theorems \ref{thm-omega} and \ref{thm-Omega}}
By (\ref{F-omega}) and (\ref{F-Omega}), we have
\begin{equation}
F_{\omega}(s, \chi)=L(s, \chi)\log L(s, \chi)-\frac{1}{2} L(s, \chi)\log L(2s, \chi^2)+L(s, \chi)G_1(s), 
\end{equation}
and 
\begin{equation}
F_{\Omega}(s, \chi)=L(s, \chi)\log L(s, \chi)+\frac{1}{2} L(s, \chi)\log L(2s, \chi^2)+L(s, \chi)G_2(s),  
\end{equation}
where $G_1(s)$ and $G_2(s)$ are absolutely convergent for $\Re (s)\geq \sigma_0=0.34$.

On $\mathcal{H}(\rho, \delta)$, by integration from $3)$ in Lemma \ref{lem-Ford-S-2.2} (see Section \ref{subsec-Lem} below), we have
\begin{equation}
\log L(s, \chi)=\log (s-\rho)+H_{\rho}(s), 
\end{equation}
with
\begin{equation}
H_{\rho}(s)=\sum_{0<|\gamma'-\gamma|\leq 1} \log(s-\rho')+O(\log |\gamma|).
\end{equation}
We write
\begin{equation}
L(s, \chi)=L'(\rho, \chi)(s-\rho)+M(s, \rho)(s-\rho)^2,
\end{equation}
where
\begin{equation}
M(s, \rho):=\int_0^1 (1-t)L''(\rho+t(s-\rho), \chi)dt.
\end{equation}
Under the ERH, we have $|M(s, \rho)|\ll |\gamma|^{\delta+\epsilon}$ on $\mathcal{H}(\rho, \delta)$. 

Now we define a function $T(x)$ as follows, for $T_{n'}\in \mathcal{T}$ satisfying $\big[e^{{\frac{5}{6}(\frac{5}{4})}^{n}}\big] \leq T_{n'}\leq \big[e^{\frac{5}{6}{(\frac{5}{4})}^{n}}\big]+1$, let $T(x)=T_{n'}$ for $e^{{(\frac{5}{4})}^n}\leq x\leq e^{{(\frac{5}{4})}^{n+1}}$. In particular, we have
$$\frac{1}{2}x^{\frac{2}{3}}\leq T(x)\leq 2x^{\frac{5}{6}} \quad (x\geq e^{\frac{5}{4}}).$$

Thus, by Lemma \ref{lem-psi-contour}, for $T=T(x)$, 
\begin{align}\label{psi-omega}
\psi_{\omega}(x, \chi)=&\sum_{|\gamma|\leq T} \frac{1}{2\pi i}\int_{\mathcal{H}(\rho, \delta)} L(s, \chi)\log L(s, \chi)\frac{x^s}{s} ds -\frac{a(\chi)}{2}\frac{1}{2\pi i}\int_{\mathcal{H}(\frac{1}{2}, \delta)} L(s, \chi)\log L(2s, \chi^2)\frac{x^s}{s} ds\nonumber\\
& +O\(x^{\frac{1}{2}-\frac{\delta}{10}}\), 
\end{align}
and 
\begin{align}\label{psi-Omega}
\psi_{\Omega}(x, \chi)=&\sum_{|\gamma|\leq T} \frac{1}{2\pi i}\int_{\mathcal{H}(\rho, \delta)} L(s, \chi)\log L(s, \chi)\frac{x^s}{s} ds +\frac{a(\chi)}{2}\frac{1}{2\pi i}\int_{\mathcal{H}(\frac{1}{2}, \delta)} L(s, \chi)\log L(2s, \chi^2)\frac{x^s}{s} ds\nonumber\\
& +O\(x^{\frac{1}{2}-\frac{\delta}{10}}\).
\end{align}

On each truncated Hankel contour $\mathcal{H}(\rho, \delta)$, 
\begin{equation}
L(s, \chi)\log L(s, \chi)=L'(\rho, \chi)(s-\rho)\log (s-\rho)+M(s, \rho)(s-\rho)^2 \log (s-\rho)+L(s, \chi)H_{\rho}(s).
\end{equation}
Since $s=\rho$ is not singularity of $L(s, \chi)H_{\rho}(s)$, and by Lemma \ref{lem-main-term}, 
\begin{align}\label{L-L-main}
\sum_{|\gamma|\leq T} &\frac{1}{2\pi i}\int_{\mathcal{H}(\rho, \delta)} L(s, \chi)\log L(s, \chi)\frac{x^s}{s} ds\nonumber\\
=&\sum_{|\gamma|\leq T} \frac{1}{2\pi i}\int_{\mathcal{H}(\rho, \delta)} \(L'(\rho, \chi)(s-\rho) \log(s-\rho)+M(s, \rho) (s-\rho)^2 \log (s-\rho)\) \frac{x^s}{s}ds\nonumber\\
=&\frac{\sqrt{x}}{\log^2 x}\sum_{|\gamma|\leq T} \frac{L'(\rho, \chi)x^{i\gamma}}{\frac{1}{2}+i\gamma}+O\(\frac{\sqrt{x}}{\log^3 x} \)+\sum_{|\gamma|\leq T} \frac{1}{2\pi i} \int_{\mathcal{H}(\rho,\delta)} M(s, \rho) (s-\rho)^2 \log (s-\rho) \frac{x^s}{s}ds.
\end{align}

For the second integral in (\ref{psi-omega}), if $\chi$ is real, $s=\frac{1}{2}$ is a pole of $L(2s, \chi^2)$. On the truncated Hankel contour $\mathcal{H}(\frac{1}{2}, \delta)$, we have
\begin{equation}\label{log-L2s}
\log L(2s, \chi^2)=-\log\(s-\frac{1}{2}\)+H_B(s), 
\end{equation}
where $H_B(s)=O(1)$ on $\mathcal{H}(\frac{1}{2}, \delta)$. And we write
\begin{eqnarray}\label{L-half-exp}
L(s, \chi)=L\(\frac{1}{2}, \chi\)+L'\(\frac{1}{2}, \chi\)\(s-\frac{1}{2}\)+M_B\(s, \frac{1}{2}\)\(s-\frac{1}{2}\)^2, 
\end{eqnarray}
where,
\begin{equation}
M_B\(s, \frac{1}{2} \):=\int_0^1 (1-t) L''\(\frac{1}{2}+t\(s-\frac{1}{2}\), \chi\) dt=O(1), \mbox{~on~} \mathcal{H}\(\frac{1}{2}, \delta \).
\end{equation}
By (\ref{log-L2s}), (\ref{L-half-exp}), and Lemma \ref{lem-main-term}, we get
\begin{align}\label{L-L2s-integral}
\frac{1}{2\pi i}&\int_{\mathcal{H}(\frac{1}{2}, \delta)} L(s, \chi)\log L(2s, \chi^2)\frac{x^s}{s} ds\nonumber\\
=&-\frac{1}{2\pi i}\int_{\mathcal{H}(\frac{1}{2}, \delta)} \(L\(\frac{1}{2}, \chi\)+L'\(\frac{1}{2}, \chi\)\(s-\frac{1}{2}\) \) \log \( s-\frac{1}{2}\) \frac{x^s}{s} ds\nonumber\\
& -\frac{1}{2\pi i}\int_{\mathcal{H}(\frac{1}{2}, \delta)} M_B\(s, \frac{1}{2}\)\(s-\frac{1}{2}\)^2 \log \( s-\frac{1}{2}\)\frac{x^s}{s} ds\nonumber\\
=& 2L\(\frac{1}{2}, \chi\) \frac{\sqrt{x}}{\log x} +\(4L\(\frac{1}{2}, \chi\)-2L'\(\frac{1}{2}, \chi\)\) \frac{\sqrt{x}}{\log^2 x}+O\( \frac{\sqrt{x}}{\log^3 x}\)\nonumber\\
& -\frac{1}{2\pi i}\int_{\mathcal{H}(\frac{1}{2}, \delta)} M_B\(s, \frac{1}{2}\)\(s-\frac{1}{2}\)^2 \log \( s-\frac{1}{2}\)\frac{x^s}{s} ds. 
\end{align}
Then, by Lemma \ref{lem-Meng-8}, since $r_0=\frac{1}{x}$, 
\begin{align}\label{L-half-error}
&\frac{1}{2\pi i}\int_{\mathcal{H}(\frac{1}{2}, \delta)} M_B\(s, \frac{1}{2}\)\(s-\frac{1}{2}\)^2 \log \( s-\frac{1}{2}\)\frac{x^s}{s} ds\nonumber\\
&\ll \left| \int_{r_0}^{\delta} M_B\(\frac{1}{2}-\sigma, \frac{1}{2}\) \sigma^2 \((\log\sigma-i\pi)-(\log\sigma+i\pi) \)\frac{x^{1/2-\sigma}}{1/2-\sigma}d\sigma \right| \nonumber\\
& \qquad + \int_0^{2\pi} r_0^2\( \log \frac{1}{r_0}+\pi\) r_0 d\alpha \nonumber\\
&\ll \frac{\sqrt{x}}{\log^3 x}+\frac{\log x}{x^3}\ll \frac{\sqrt{x}}{\log^3 x}.
\end{align}
By (\ref{L-L2s-integral}) and (\ref{L-half-error}), we deduce that
\begin{align}\label{L-L2s-formula}
\frac{1}{2\pi i}&\int_{\mathcal{H}(\frac{1}{2}, \delta)} L(s, \chi)\log L(2s, \chi^2)\frac{x^s}{s} ds\nonumber\\
=& 2L\(\frac{1}{2}, \chi\) \frac{\sqrt{x}}{\log x} +\(4L\(\frac{1}{2}, \chi\)-2L'\(\frac{1}{2}, \chi\)\) \frac{\sqrt{x}}{\log^2 x}+O\( \frac{\sqrt{x}}{\log^3 x}\). 
\end{align}

Thus, by (\ref{psi-omega}), (\ref{psi-Omega}), (\ref{L-L-main}), and (\ref{L-L2s-formula}), we have
\begin{align}\label{psi-omega-formula}
\psi_{\omega}(x, \chi)=&-a(\chi)\left\{L\(\frac{1}{2}, \chi\) \frac{\sqrt{x}}{\log x} +\(2L\(\frac{1}{2}, \chi\)-L'\(\frac{1}{2}, \chi\)\) \frac{\sqrt{x}}{\log^2 x}\right\}\nonumber\\
&+\frac{\sqrt{x}}{\log^2 x}\sum_{|\gamma|\leq T} \frac{L'(\rho, \chi)x^{i\gamma}}{\frac{1}{2}+i\gamma}+\sum_{|\gamma|\leq T} \frac{1}{2\pi i} \int_{\mathcal{H}(\rho,\delta)} M(s, \rho) (s-\rho)^2 \log (s-\rho) \frac{x^s}{s}ds\nonumber\\
&+O\( \frac{\sqrt{x}}{\log^3 x}\), 
\end{align}
and 
\begin{align}\label{psi-Omega-formula}
\psi_{\Omega}(x, \chi)=&a(\chi)\left\{L\(\frac{1}{2}, \chi\) \frac{\sqrt{x}}{\log x} +\(2L\(\frac{1}{2}, \chi\)-L'\(\frac{1}{2}, \chi\)\) \frac{\sqrt{x}}{\log^2 x}\right\}\nonumber\\
&+\frac{\sqrt{x}}{\log^2 x}\sum_{|\gamma|\leq T} \frac{L'(\rho, \chi)x^{i\gamma}}{\frac{1}{2}+i\gamma}+\sum_{|\gamma|\leq T} \frac{1}{2\pi i} \int_{\mathcal{H}(\rho,\delta)} M(s, \rho) (s-\rho)^2 \log (s-\rho) \frac{x^s}{s}ds\nonumber\\
&+O\( \frac{\sqrt{x}}{\log^3 x}\).
\end{align}

In the following, for $T=T(x)$, define
\begin{equation}
\Sigma(x, \chi):=\frac{\log^2 x}{\sqrt{x}}\sum_{|\gamma|\leq T} \int_{\mathcal{H}(\rho,\delta)} M(s, \rho) (s-\rho)^2 \log (s-\rho) \frac{x^s}{s}ds.
\end{equation}
For fixed large $T_0$, let
\begin{equation}
S(x, T_0; \chi):=\sum_{|\gamma|\leq T(x)}\frac{L'(\rho, \chi)x^{i\gamma}}{\frac{1}{2}+i\gamma}-\sum_{|\gamma|\leq T_0}\frac{L'(\rho, \chi)x^{i\gamma}}{\frac{1}{2}+i\gamma}.
\end{equation}

Then we have the following results. We give their proof in later subsection.
\begin{lem}\label{lem-main-error}
\begin{equation}
\int_2^Y |\Sigma(e^y, \chi)|^2 dy=O(1).
\end{equation}	
\end{lem}
The following lemma is similar to Lemma 11 in \cite{Meng}.
\begin{lem}\label{lem-Meng-13}
	For any $\epsilon>0$, 
	\begin{equation}
	\int_2^Y |S(e^y, T_0; \chi)|^2 dy\ll \frac{Y}{T_0^{1-\epsilon}} + \frac{\log Y}{T_0^{1-\epsilon}}+T_0^{\epsilon}.
	\end{equation}
\end{lem}

 By (\ref{psi-omega-formula}), (\ref{psi-Omega-formula}), and Lemmas \ref{lem-main-error} and \ref{lem-Meng-13},  Theorems \ref{thm-omega} and \ref{thm-Omega} follow. \qed

\vspace{1em}

\noindent\textbf{Remark.} The result of Lemma \ref{lem-main-error} seems good. Actually, we have
\begin{equation}
\int_2^X |\Sigma(x, \chi)|^2 dx=o(X).
\end{equation}
However, Lemma \ref{lem-Meng-13} prevents us from proving a better result for our theorem.

\subsection{Proof of Lemmas \ref{lem-main-error} and \ref{lem-Meng-13}}\label{subsec-Lem}
First, we need the following lemmas on the zeros of Dirichlet $L$-functions.

Let $\gamma$ be the imaginary part of a zero of $L(s, \chi)$ in the critical strip.
\begin{lem}[\cite{Ford-S}, Lemma 2.2] \label{lem-Ford-S-2.2}
	Let $\chi$ be a Dirichlet character modulo $q$. Let $N(T, \chi)$ denote the number of zeros of $L(s, \chi)$ with $0<\Re(s)<1$ and $|\Im(s)|<T$. Then
	
	1) $N(T, \chi)=O(T\log (qT))$ for $T\geq 1$.
	
	2) $N(T, \chi)-N(T-1, \chi)=O(\log (qT))$ for $T\geq 1$.
	
	3) Uniformly for $s=\sigma+it$ and $\sigma\geq -1$,
	\begin{equation}\label{F-S-lem-formu}
	\frac{L'(s, \chi)}{L(s, \chi)}=\sum_{|\gamma-t|<1} \frac{1}{s-\rho} +O(\log q(|t|+2)).
	\end{equation}
\end{lem}

Similar to Lemma 2.4 in \cite{Ford-S}, we have
\begin{lem} \label{lem-Ford-S-2.4}
	Assume $L(\frac{1}{2}, \chi)\neq 0$. For $A\geq 0$ and $0\leq \delta<\frac{1}{2}$,
	\begin{equation*}
	\sum_{\substack{|\gamma_1|, |\gamma_2|\geq A\\ |\gamma_1-\gamma_2|\geq 1}} \frac{1}{|\gamma_1|^{1-\delta}|\gamma_2|^{1-\delta}|\gamma_1-\gamma_2|}\ll \frac{\log(A+2)}{(A+1)^{1-2\delta}}.
	\end{equation*}
\end{lem}
\noindent\textit{\textbf{Proof of Lemma \ref{lem-Ford-S-2.4}.}} We only need to estimate the case $|\gamma_2|\geq |\gamma_1|$. By Lemma \ref{lem-Ford-S-2.2} and partial summation, the sum in the lemma is
\begin{align}
&\ll \sum_{|\gamma_2|\geq A} \frac{1}{|\gamma_2|^{1-\delta}}\(\frac{1}{|\gamma_2|} \sum_{|\gamma_1|<\frac{|\gamma_2|}{2}}\frac{1}{|\gamma_1|^{1-\delta}}+\frac{1}{|\gamma_2|^{1-\delta}}\sum_{\substack{\frac{|\gamma_2|}{2}\leq |\gamma_1|\leq |\gamma_2|\\ |\gamma_2-\gamma_1|\geq 1}} \frac{1}{|\gamma_2-\gamma_1|} \)\nonumber\\
&\ll \sum_{|\gamma_2|\geq A} \frac{1}{|\gamma_2|^{1-\delta}} \frac{\log (|\gamma_2|+2)}{|\gamma_2|^{1-\delta}}\ll \frac{\log (A+2)}{(A+1)^{1-2\delta}}.
\end{align}
We get the conclusion of this lemma. \qed

\bigskip

In the following, we use the above lemmas to prove Lemmas \ref{lem-main-error} and \ref{lem-Meng-13}.

\noindent \textit{\textbf{Proof of Lemma \ref{lem-main-error}.}}
We rewrite $\Sigma(x, \chi)$ as 
\begin{equation}
\log^2 x\sum_{|\gamma|\leq T} \int_{\mathcal{H}(\rho, \delta)}x^{i\gamma} M(s, \rho) (s-\rho)^2 \log (s-\rho)\frac{x^{s-\rho}}{s}ds.
\end{equation}

Let $\Gamma_{\rho}$ represent the circle in the Hankel contour $\mathcal{H}(\rho, \delta)$. 
Denote
\begin{align}\label{ExRho-h-r}
E(x; \rho):=&\int_{\mathcal{H}(\rho, \delta)}M(s, \rho) (s-\rho)^2 \log (s-\rho)\frac{x^{s-\rho}}{s}ds\nonumber\\
=& \int_{r_{\rho}}^{\delta} M(1/2-\sigma+i\gamma, \rho) \sigma^2 \((\log\sigma-i\pi)-(\log\sigma+i\pi) \)\frac{x^{-\sigma}}{1/2-\sigma+i\gamma}d\sigma\nonumber\\
&+\int_{\Gamma_{\rho}}M(s, \rho) (s-\rho)^2 \log (s-\rho)\frac{x^{s-\rho}}{s}ds\nonumber\\
=:&  E_h(x; \rho)+E_r(x; \rho).
\end{align}

Since we assume $r_{\rho}\ll \frac{1}{x}$, and $T(x)\ll x^{\frac{5}{6}}$, 
\begin{equation}\label{E-r-square}
\left|\sum_{|\gamma|\leq T} x^{i\gamma}E_r(x; \rho) \right|^2\ll \sum_{|\gamma|\leq T}|E_r(x; \rho)|\ll \sum_{|\gamma|\leq T}\frac{1}{|\gamma|^{1-\delta-\epsilon}} \frac{1}{x^2} r_{\rho} \(\log \frac{1}{r_\rho}+\pi\)\ll \frac{1}{x^{2-\epsilon}}.
\end{equation}

Next, we take care of the following sum
\begin{align}
\left|\sum_{|\gamma|\leq T} x^{i\gamma}E_h(x; \rho)\right|^2&=\(\sum_{\substack{|\gamma_1-\gamma_2|\leq 1\\|\gamma_1|, |\gamma_2|\leq T}}+\sum_{\substack{|\gamma_1-\gamma_2|> 1\\|\gamma_1|, |\gamma_2|\leq T}}\)  x^{i(\gamma_1-\gamma_2)} E_h(x; \rho_1) E_h(x; \overline{\rho}_2)\nonumber\\
&=:\Sigma_1(x)+\Sigma_2(x).
\end{align}
By (\ref{ExRho-h-r}) and using change of variable, 
\begin{equation}
|E_h(x; \rho)|\ll \frac{1}{|\gamma|^{1-\delta-\epsilon}}\int_0^{\delta} \sigma^2 x^{-\sigma}d\sigma\ll \frac{1}{|\gamma|^{1-\delta-\epsilon}} \frac{1}{\log^3 x}
\end{equation}
Hence, 
\begin{equation}\label{dist-error-sum1}
\int_2^Y y^4 \Sigma_1(e^y) dy \ll \int_2^Y \frac{1}{y^2}\sum_{\gamma} \frac{1}{|\gamma|^{2-2\delta-\epsilon}}dy=O(1). 
\end{equation}

For $\Sigma_2(x)$,
\begin{equation}
\Sigma_2(x)=\sum_{\substack{|\gamma_1-\gamma_2|>1\\ |\gamma_1|, |\gamma_2|\leq T}}x^{i(\gamma_1-\gamma_2)} E_h(x; \rho_1)  E_h(x; \overline{\rho}_2).
\end{equation}
For $e^{{(\frac{5}{4})}^l}\leq x\leq e^{{(\frac{5}{4})}^{l+1}}$, $T=T(x)=T_{l'}$ is a constant, so we define
\begin{equation}\label{pf-distri-Err-lem-J-J(x)}
J(x):= \sum_{\substack{|\gamma_1-\gamma_2|>1\\|\gamma_1|, |\gamma_2|\leq T_{l'}}} x^{i(\gamma_1-\gamma_2)} \int_{r_{\rho_1}}^{\delta} \int_{r_{\overline{\rho}_2}}^{\delta} R_{\rho_1}(\sigma_1; x) R_{\overline{\rho}_2}(\sigma_2; x) \frac{d\sigma_1 d\sigma_2}{i(\gamma_1-\gamma_2)-(\sigma_1+\sigma_2)},
\end{equation}
where
$ R_{\rho}(\sigma; x)=M(1/2-\sigma+i\gamma, \rho) \sigma^2 \((\log\sigma-i\pi)-(\log\sigma+i\pi) \)\frac{x^{-\sigma}}{1/2-\sigma+i\gamma}$.
Thus,
\begin{equation}\label{pf-distri-Err-lem-J-sum-to-J}
\int_{e^{{(\frac{5}{4})}^l}}^{e^{{(\frac{5}{4})}^{l+1}}} \sum_{\substack{|\gamma_1-\gamma_2|>1\\|\gamma_1|, |\gamma_2|\leq T_{l'}}}  x^{i(\gamma_1-\gamma_2)} E_h(x; \rho_1) E_h(x;\overline{\rho}_2)\frac{dx}{x}=J(e^{{(\frac{5}{4})}^{l+1}})-J(e^{{(\frac{5}{4})}^l}).
\end{equation}
By Lemma \ref{lem-Ford-S-2.4}, 
\begin{equation}
|J(x)|\ll \sum_{|\gamma_1-\gamma_2|>1} \frac{1}{|\gamma_1|^{1-\delta-\epsilon}|\gamma_2|^{1-\delta-\epsilon}|\gamma_1-\gamma_2|} \frac{1}{\log^6 x}\ll \frac{1}{\log^6 x}.
\end{equation}
Hence, for any positive integer $l$, 
\begin{equation}
\int_{(\frac{5}{4})^l}^{(\frac{5}{4})^{l+1}} \Sigma_2(e^y) dy\ll \frac{1}{{(5/4)}^{6l}}.
\end{equation}

Therefore, 
\begin{equation}\label{dist-error-sum2}
\int_2^Y y^4 \Sigma_2(e^y) dy\ll \sum_{l\leq \frac{\log Y}{\log(5/4)}+1} \(\frac{5}{4}\)^{4l} \int_{{(\frac{5}{4})}^l}^{{(\frac{5}{4})}^{l+1}} \Sigma_2(e^y) dy\ll \sum_{l\leq \frac{\log Y}{\log(5/4)}+1} \frac{1}{(5/4)^{2l}}=O(1).
\end{equation}

Combining (\ref{E-r-square}), (\ref{dist-error-sum1}), and (\ref{dist-error-sum2}), we get the desired result. \qed

\vspace{1em}

\noindent\textit{\textbf{Proof of Lemma \ref{lem-Meng-13}.}} For fixed $T_0$, let $X_0$ be the largest $x$ such that $T=T(x)\leq T_0$. Since $x^{\frac{2}{3}}\ll T(x)\ll x^{\frac{5}{6}}$, $\log X_0 \asymp \log T_0$. Under the ERH, $|L'(\rho, \chi)|\ll |\gamma|^{\epsilon}$, by Lemmas \ref{lem-Ford-S-2.2} and \ref{lem-Ford-S-2.4},
\begin{align*}
	&\int_2^Y | S(e^y, T_0; \chi)|^2 dy\leq \int_2^{\log X_0} \left| \sum_{|\gamma|\leq T_0} \frac{1}{|\gamma|^{1-\epsilon}} \right|^2 dy +\int_{\log X_0}^Y \left| \sum_{T_0\leq |\gamma|\leq T(e^y)} \frac{L'(\rho, \chi)e^{i\gamma y}}{\frac{1}{2}+i\gamma} \right|^2 dy \nonumber\\
	& \ll T_0^{\epsilon} + \sum_{\frac{\log\log X_0}{\log(5/4)}\leq l\leq \frac{\log Y}{\log(5/4)}+1} \int_{(\frac{5}{4})^l}^{(\frac{5}{4})^{l+1}} \( \sum_{\substack{|\gamma_1-\gamma_2|\leq 1\\ T_0\leq |\gamma_1|, |\gamma_2|\leq T_{l'}}} +\sum_{\substack{|\gamma_1-\gamma_2|> 1\\ T_0\leq |\gamma_1|, |\gamma_2|\leq T_{l'}}} \) \frac{L'(\rho_1, \chi) L'(\bar{\rho}_2, \chi)e^{i(\gamma_1-\gamma_2)y}}{(\frac{1}{2}+i\gamma_1)(\frac{1}{2}-i\gamma_2)} dy \nonumber\\
	&\ll T_0^{\epsilon}+ \sum_{\frac{\log\log X_0}{\log(5/4)}\leq l\leq \frac{\log Y}{\log(5/4)}+1} \( \(\frac{5}{4}\)^l\sum_{|\gamma|\geq T_0} \frac{1}{|\gamma|^{2-\epsilon}}+ \sum_{|\gamma_1|, |\gamma_2| \geq T_0} \frac{1}{|\gamma_1|^{1-\epsilon}|\gamma_2|^{1-\epsilon}|\gamma_1-\gamma_2|}\)\nonumber\\
	& \ll \frac{Y}{T_0^{1-\epsilon}} + \frac{\log Y}{T_0^{1-\epsilon}}+T_0^{\epsilon}.
\end{align*}
This completes the proof of this lemma. \qed

\section{General Theorem in the Selberg class}\label{Sec-Selberg}
Selberg class $\mathcal{S}$ consists a class of Dirichlet series introduced by Selberg \cite{Selberg}. A Dirichlet series
\begin{equation}
F(s)=\sum_{n=1}^{\infty} \frac{a(n)}{n^s} \quad (\sigma>1),
\end{equation}
is in $\mathcal{S}$ provided it satisfies the following hypotheses. 

1) Ramanujan Hypothesis: $a(n)\ll n^{\epsilon}$ for any $\epsilon>0$;

2) Analyticity: $(s-1)^m F(s)$ is an entire function  of finite order for some nonnegative integer $m$; 

3)Functional Equation: there exist positive real numbers $Q$, $\lambda_1, \cdots, \lambda_r$, and complex numbers $w, \mu_1, \cdots, \mu_r$ with $\Re \mu_j\geq 0$ and $|w|=1$ such that
$$ \Lambda_F(s)=w \overline{\Lambda_F(1-\bar{s})},$$
where
$$\Lambda_F(s)=F(s)Q^s \prod_{j=1}^r \Gamma(\lambda_j s+\mu_j);$$

4) Euler Product: for prime powers $p^k$ there exist complex numbers $b(p^k)$ satisfying 
$$\log F(s)=\sum_p \sum_{m=1}^{\infty} \frac{b(p^{m})}{p^{ms}}$$
and $b(p^m)\ll p^{m\theta}$ for some $\theta<\frac{1}{2}$. 

Let $d_F$ be the \textit{degree} of $F\in\mathcal{S}$ defined by 
$$d_F=2\sum_{j=1}^r \lambda_j.$$

We assume the corresponding Hypotheses, 1) Grand Riemann Hypothesis (GRH): if $F\in \mathcal{S}$, then $F(s)\neq 0$ for $\sigma>\frac{1}{2}$; 
 and the Lindel\"{o}f Hypothesis ($\rm LH_S$): $\forall \epsilon>0$, $|F(s)|\ll (|t|+2)^{\epsilon}$ for $\frac{1}{2}\leq \sigma\leq 1$.  
 
 If the Dirichlet series $\sum_{n=1}^{\infty} \frac{f(n)}{n^s}$ can be written as linear combinations of the form $\log^k F(s)$ for $F\in\mathcal{S}$ and $F(s)$ has no pole, then we can prove the existence of a limiting distribution related to the summatory function
\begin{equation}
S_f(x):=\sum_{n\leq x} f(n).
\end{equation}

It suffices to consider
\begin{equation}
\frac{1}{2\pi i}\int_{a-i\infty}^{a+i\infty} \log^k F(s) \frac{x^s}{s} ds.
\end{equation}
Suppose
\begin{equation}
\sum_{n=1}^{\infty} \frac{s(n)}{n^s}=\log^k F(s),
\end{equation}
let
\begin{equation}
S(x):=\sum_{n\leq x} s(n).
\end{equation}

We have the following theorem.
\begin{theorem}\label{thm-SC} Under the  assumptions GRH and $ LH_S$, 
	\begin{equation}
	S(x)=\frac{(-1)^k k\sqrt{x}}{\log x}(\log\log x)^{k-1} \left\{ \sum_{|\gamma|\leq T_0}\frac{ m^k(\gamma) x^{i\gamma}}{\frac{1}{2}+i\gamma}+\Sigma(x, T_0)\right\}, 
	\end{equation}
\end{theorem}
where
\begin{equation}
\limsup_{Y\rightarrow\infty}\frac{1}{Y}\int_1^Y \left| \Sigma(e^y, T_0) \right|^2 dy\ll\frac{(\log T_0)^{2k+2}}{T_0}.
\end{equation}

\subsection{Sketch of the proof of Theorem \ref{thm-SC}}
We have corresponding results for Selberg class. 
\begin{lem}[\cite{Kac-Pere}, (2.4), or \cite{Pank-Steding}, (4)]\label{lem-SC-zero-density}
	Let $F(s)$ be in Selberg class. Let $N_F(T)$ denote the number of zeros of $F(s)$ with $0<\Re (s)<1$ and $0<\Im (s)<T$. Then
	\begin{equation}
	N_F(T)=\frac{d_F}{2\pi} T\log T+c_F T+O(\log T),
	\end{equation}
	where $d_F$ is the degree of $F(s)$ and $c_F$ is a positive constant.
\end{lem}

\begin{lem}[\cite{Pank-Steding}, Lemma 2.1]\label{lem-SC-log-L}
	Let $F(s)$ be in the Selberg class and denote the nontrivial zeros of $F(s)$ by $\rho=\beta+i\gamma$. Then, for $-\frac{5}{2}\leq \sigma\leq \frac{7}{2}$, 
	\begin{equation}
	\log F(s)=\sum_{|t-\gamma|\leq 1} \log(s-\rho)+O(\log |t|),
	\end{equation}
	where $-\pi<\Im \log(s-\rho)\leq \pi$ for any $t$ not equal to an ordinate of a zero of $F(s)$.
\end{lem}

By Lemma \ref{lem-SC-zero-density}, for any integer $n\geq 0$, there exists $n\leq T_n\leq n+1$ such that the distance of $\sigma+iT_n$ to the nearest zero is $\gg_{d_F} \frac{1}{\log n}$ for all $-1<\sigma<2$. Then, by Lemma \ref{lem-SC-log-L}, for $-1<\sigma<2$, 
\begin{equation}
|\log F(\sigma+iT_n)|\ll_{d_F} \log n\log\log n\ll_{d_F} \log^2 T_n.
\end{equation}
Denote $\mathcal{T}:=\{T_n \}_{n=0}^{\infty}$. 

Then, using the same type of integration contour as Figure \ref{contour}, we have
\begin{lem}\label{lem-SC-Perron}
	Assume the GRH, $LH_S$, and $F(\frac{1}{2})\neq 0$. Let $\rho$ be nontrivial zero of $F(s)$. Then for $T\in \mathcal{T}$, 
	\begin{equation}
	S(x)=\sum_{|\gamma|\leq T} \frac{1}{2\pi i} \int_{\mathcal{H}(\rho, \delta)} \log^k F(s)\frac{x^s}{s}ds +O_k\(\frac{x\log x}{T}+\frac{x(\log T)^{2k}}{T}+x^{\frac{1}{2}-\delta}(\log T)^{k+1} \).
	\end{equation}	
\end{lem}

By Lemma \ref{lem-SC-log-L}, on $\mathcal{H}(\rho, \delta)$, we have
\begin{align}\label{SC-Hankel}
&\sum_{|\gamma|\leq T} \frac{1}{2\pi i} \int_{\mathcal{H}(\rho, \delta)} \log^k F(s)\frac{x^s}{s}ds\nonumber\\
& =\sum_{|\gamma|\leq T} \frac{1}{2\pi i} \int_{\mathcal{H}(\rho, \delta)} (m(\gamma)\log (s-\rho))^k \frac{x^s}{s}ds\nonumber\\
&\quad + \sum_{|\gamma|\leq T} \frac{1}{2\pi i} \int_{\mathcal{H}(\rho, \delta)} \sum_{j=1}^k {k \choose j}(m(\gamma)\log (s-\rho))^{k-j} (H_{\rho}(s))^j \frac{x^s}{s}ds,
\end{align}
where $m(\gamma)$ is the multiplicity of the zero $\rho$ and $m(\gamma)\ll \log |\gamma|$ by Lemma \ref{lem-SC-zero-density}, and 
\begin{equation}
H_{\rho}(s)=\sum_{0<|\gamma'-\gamma|\leq 1} m(\gamma')\log (s-\rho')+O(\log |\gamma|) \quad \mbox{on~} \mathcal{H}(\rho, \delta).
\end{equation}

For the first integral in (\ref{SC-Hankel}), by Lemma \ref{lem-main-term} and Lemma \ref{lem-SC-zero-density}, we get
\begin{align}\label{SC-main-term}
&\sum_{|\gamma|\leq T} \frac{1}{2\pi i} \int_{\mathcal{H}(\rho, \delta)} (m(\gamma)\log (s-\rho))^k \frac{x^s}{s}ds\nonumber\\
&=\frac{(-1)^k \sqrt{x}}{\log x}\(\sum_{j=1}^k {k\choose j} (\log\log x)^{k-j} \frac{1}{\Gamma_j(0)} \) \sum_{|\gamma|\leq T} \frac{m^k(\gamma) x^{i\gamma}}{\frac{1}{2}+i\gamma}+O_k\(\frac{\sqrt{x}(\log\log x)^{k-1}}{\log^2 x}\).
\end{align} 

For the second integral in (\ref{SC-Hankel}), with the same $T(x)$ as before, using the same proof as Lemma 14 in \cite{Meng}, we have the following lemma. 
\begin{lem}\label{lem-SC-mean-square}
	Let $\rho$ be a zero of $F(s)$. Assume the function $g(s)\ll \( \log |\gamma|\)^c$ on $\mathcal{H}(\rho, \delta)$ for some constant $c\geq 0$, and
	\begin{equation}\label{lem-distri-sum-func-H}
	H_{\rho}(s)=\sum_{0<|\gamma'-\gamma|\leq 1} m(\gamma')\log (s-\rho') +O\(\log |\gamma| \)\quad{\rm on~} \mathcal{H}(\rho,\delta).
	\end{equation}
	For any integers $l, n\geq 0$, denote
	\begin{equation*}
	E(x; \rho):=\int_{\mathcal{H}(\rho,\delta)} \(m(\gamma)\log(s-\rho)\)^{l} \(H_{\rho}(s)\)^n g(s) \frac{x^{s-\rho}}{s}ds.
	\end{equation*}
	Then, for $T=T(x)$, we have
	\begin{equation*}
	\int_{2}^Y \left|y \sum_{|\gamma|\leq T(e^y)} e^{i\gamma y} E(e^y; \rho) \right|^2 dy=o\( Y(\log Y)^{2l+2n-2} \).
	\end{equation*}
\end{lem}

Denote
\begin{equation}
S_1(x):=(-1)^k\(\sum_{j=2}^k {k\choose j} (\log\log x)^{k-j} \frac{1}{\Gamma_j(0)} \) \sum_{|\gamma|\leq T} \frac{m^k(\gamma) x^{i\gamma}}{\frac{1}{2}+i\gamma},
\end{equation}
and 
\begin{equation}
S_2(x; T_0):=\sum_{|\gamma|\leq T(x)}\frac{m^k(\gamma)x^{i\gamma}}{\frac{1}{2}+i\gamma}-\sum_{|\gamma|\leq T_0}\frac{m^k(\gamma)x^{i\gamma}}{\frac{1}{2}+i\gamma}.
\end{equation}
Similar to Lemma 11 in \cite{Meng}, we have
\begin{lem}\label{lem-SC-T0-square}
	\begin{equation}
	\int_2^Y |S_1(e^y)|^2 dy=o(Y(\log Y)^{2k-2}), 
	\end{equation}
	and for fixed large $T_0$, 
	\begin{equation}
	\int_2^Y |S_2(e^y, T_0)|^2 dy\ll Y \frac{(\log T_0)^{2k+2}}{T_0}+\log Y \frac{(\log T_0)^{2k+3}}{T_0}+ (\log T_0)^{2k+5}. 
	\end{equation}
\end{lem}

Combining (\ref{SC-Hankel}) and (\ref{SC-main-term}) with Lemmas \ref{lem-SC-mean-square} and \ref{lem-SC-T0-square}, Theorem \ref{thm-SC} follows. \qed

\bigskip
\textbf{Acknowledgement.} I would like to thank Professor Greg Martin for proposing me this project and kindly providing me related numerical data.

{\footnotesize
\noindent Centre de Recherches Math\'{e}matiques,
Universit\'{e} de Montr\'{e}al,
Montr\'{e}al, QC, H3C 3J7 Canada

\noindent \textit{E-mail}: meng@crm.umontreal.ca

\end{document}